\documentclass[12pt,a4paper]{article}
\usepackage[cp1251]{inputenc}
\usepackage{amsfonts, amssymb, amsthm, amsmath}
\usepackage{mathtext}
\usepackage{longtable}
\usepackage[russian]{babel}

\theoremstyle{remark}

\renewcommand{\le}{\leqslant}
\renewcommand{\ge}{\geqslant}

\renewcommand{\phi}{\varphi}

\newcommand{\eqd}{\stackrel{d}{=}}

\title{{On convergence of the distributions of random sequences
with independent random indexes to variance-mean
mixtures}\thanks{Research supported by the Russian Scientific
Foundation (project 14-11-00364).}}
\author{V. Yu. Korolev\thanks{Faculty of
Computational Mathematics and Cybernetics, Lomonosov Moscow State
University; Institute of Informatics Problems, FRC IC RAS;
victoryukorolev@yandex.ru}, A. I. Zeifman\thanks{Vologda State
University; Institute of Informatics Problems, FRC IC RAS; ISEDT RAS;
a$\_$zeifman@mail.ru}}

\date{}

\begin{document}

\maketitle

{\small

{\bf Abstract.} We prove a version of a general transfer theorem for
random sequences with independent random indexes in the double array
limit setting under relaxed conditions. We also prove its partial
inverse providing the necessary and sufficient conditions for the
convergence of randomly indexed random sequences. Special attention
is paid to the case where the elements of the basic double array are
formed as cumulative sums of independent not necessarily identically
distributed random variables. Using simple moment-type conditions we
prove the theorem on convergence of the distributions of such sums
to normal variance-mean mixtures.

\smallskip


}

\section{Introduction}

Random sequences with independent random indexes play an important
role in modeling real processes in many fields. Most popular
examples of the application of these models usually deal with
insurance and reliability theory \cite{Kalashnikov1997,
BeningKorolev2002}, financial mathematics and queuing theory
\cite{BeningKorolev2002, GnedenkoKorolev1996}, chaotic processes in
plasma physics \cite{KorolevSkvortsova2006} where random sums are
the principal mathematical models. More general randomly indexed
random sequences arrive in the statistics of samples with random
sizes. Indeed, very often the data to be analyzed is collected or
registered during a certain period of time and the flow of
informative events producing the observations forms a random point
process, so that the number of available observations is unknown
till the end of the process of their registration and also must be
treated as a (random) observation.

The randomness of indexes usually leads to that the limit
distributions for the corresponding randomly indexed random
sequences are heavy-tailed even in the situations where the
distributions of non-randomly indexed random sequences are
asymptotically normal see, e. g., \cite{GnedenkoKorolev1996,
BeningKorolev2003}.

The literature on random sequences with random indexes is extensive,
see, e. g., the references above and the references therein. The
mathematical theory of random sequences with random indexes and, in
particular, random sums, is well-developed. However, there still
remain some unsolved problems. For example, the necessary and
sufficient conditions for the convergence of the distributions of
random sums to normal variance-mean mixtures were found only
recently for the case of identically distributed summands (see
\cite{Korolev2013, GK2013}). The case of random sums of {\it
non-identically} distributed random summands and, moreover, more
general statistics constructed from samples with random sizes has
not been considered yet. At the same time, normal variance-mean
mixtures are widely used as mathematical models of statistical
regularities in many fields. In particular, in 1977--78
O.\,Barndorff-Nielsen \cite{BN1977}, \cite{BN1978} introduced the
class of {\it generalized hyperbolic distributions} as a class of
special variance-mean mixtures of normal laws in which the mixing is
carried out with respect to one parameter since the location and
scale parameters of the mixed normal distribution are directly
linked. The range of applications of generalized hyperbolic
distributions varies from the theory of turbulence or particle size
description to financial mathematics, see \cite{BN1982}.

The paper is organized as follows. Basic notation is introduced in
Section 2. Here an auxiliary result on the asymptotic rapprochement
of the distributions of randomly indexed random sequences with
special scale-location mixtures is proved. In Section 3 of the
present paper we prove a a version of a general transfer theorem for
random sequences with independent random indexes in the double array
limit setting under relaxed conditions. We also prove its partial
inverse providing the necessary and sufficient conditions for the
convergence of randomly indexed random indexes. Following the lines
of \cite{Korolev1993}, we first formulate a general result improving
some results of \cite{Korolev1993,BeningKorolev2002} by removing
some superfluous assumptions and relaxing some conditions. The
special attention is paid to the case where the elements of the
basic double array are formed as cumulative sums of independent not
necessarily identically distributed random variables. This case is
considered in Section 4. To prove our results, we use simply
tractable moment-type conditions which can be easily interpreted
unlike general conditions providing the weak convergence of random
sums of non-identically distributed summands in \cite{Szasz1972,
KruglovKorolev1990} and \cite{KruglovZhangBo2001}. In Section 5 we
prove the theorem on convergence of the distributions of such sums
to normal variance-mean mixtures. As a simple corollary of this
result we can obtain some results of the recent paper
\cite{Toda2011}. That paper demonstrates that there is still a
strong interest to geometric sums of non-identically distributed
summands and to the application of the skew Laplace distribution
which is a normal variance-mean mixture under exponential mixing
distribution \cite{KorolevSokolov2012}.

\section{Notation. Auxiliary results}

Assume that all the random variables considered in this paper are
defined on one and the same probability space
$(\Omega,\,\mathfrak{F},\,{\sf P})$. In what follows the symbols
$\eqd$ and $\Longrightarrow$ will denote coincidence of
distributions and weak convergence (convergence in distribution). A
family $\{X_j\}_{j\in\mathbb{N}}$ of random variables is said to be
{\it weakly relatively compact}, if each sequence of its elements
contains a weakly convergent subsequence. In the finite-dimensional
case the weak relative compactness of a family
$\{X_j\}_{j\in\mathbb{N}}$ is equivalent to its {\it tightness}:
$$
\lim_{R\to\infty}\sup_{n\in\mathbb{N}}{\sf P}(|X_n|>R)=0
$$
(see, e. g., \cite{Loeve1977}).

Let $\{S_{n,k}\}$, $n,k\in\mathbb{N},$ be a double array of random
variables. For $n,k\in\mathbb{N}$ let $a_{n,k}$ and $b_{n,k}$ be
real numbers such that $b_{n,k}>0$. The purpose of the constants
$a_{n,k}$ and $b_{n,k}$ is to provide weak relative compactness of
the family of the random variables
$$
\Big\{Y_{n,k}\equiv\frac{S_{n,k}-a_{n,k}}{b_{n,k}}\Big\}_{n,k\in\mathbb{N}}
$$
in the cases where it is required.

Consider a family $\{N_n\}_{n\in\mathbb{N}}$ of nonnegative integer
random variables such that for each $n,k\in\mathbb{N}$ the random
variables $N_n$ and $S_{n,k}$ are independent. Especially note that
we do not assume the row-wise independence of $\{S_{n,k}\}_{k\ge1}$.
Let $c_n$ and $d_n$ be real numbers, $n\in\mathbb{N}$, such that
$d_n>0$. Our aim is to study the asymptotic behavior of the random
variables
$$
Z_n\equiv\frac{S_{n,N_n}-c_n}{d_n}
$$
as $n\to\infty$ and find rather simple conditions under which the
limit laws for $Z_n$ have the form of normal variance-mean mixtures.
In order to do so we first formulate a somewhat more general result
following the lines of \cite{Korolev1993}, removing superfluous
assumptions, relaxing the conditions and generalizing some of the
results of that paper.

The characteristic functions of the random variables $Y_{n,k}$ and
$Z_n$ will be denoted $h_{n,k}(t)$ and $f_n(t)$, respectively,
$t\in\mathbb{R}$.

Let $Y$ be a random variable whose distribution function and
characteristic function will be denoted $H(x)$ and $h(t)$,
respectively, $x,t\in\mathbb{R}$. Introduce the random variables
$$
U_n=\frac{b_{n,N_n}}{d_n},\ \ \ V_n=\frac{a_{n,N_n}-c_n}{d_n}.
$$
Introduce the function
$$
g_n(t)\equiv {\sf
E}h(tU_n)\exp\{itV_n\}=\sum\nolimits_{k=1}^{\infty}{\sf
P}(N_n=k)\exp\Big\{it\frac{a_{n,k}-c_n}{d_n}\Big\}h\Big(\frac{tb_{n,k}}{d_n}\Big),\
\ \ t\in\mathbb{R}.
$$
It can be easily seen that $g_n(t)$ is the characteristic function
of the random variable $Y\cdot U_n+V_n$ where the random variable
$Y$ is independent of the pair $(U_n,\,V_n)$. Therefore, the
distribution function $G_n(x)$ corresponding to the characteristic
function $g_n(t)$ is the scale-location mixture of the distribution
function $H(x)$:
$$
G_n(x)={\sf E}H\Big(\frac{x-V_n}{U_n}\Big),\ \ \
x\in\mathbb{R}.\eqno(1)
$$

In the double-array limit setting considered in this paper, to
obtain non-trivial limit laws for $Z_n$ we require the following
additional {\it coherency condition}: for any $T\in(0,\infty)$
$$
\lim_{n\to\infty}{\sf E}\sup_{|t|\le
T}\big|h_{n,N_n}(t)-h(t)\big|=0.\eqno(2)
$$
To clarify the sense of the coherency condition, note that if we had
usual row-wise convergence of $Y_{n,k}$ to $Y$, then for any
$n\in\mathbb{N}$ and $T\in[0,\infty)$
$$
\lim_{k\to\infty}\sup_{|t|\le T}|h_{n,k}(t)-h(t)|=0.\eqno(3)
$$
So we can say that coherency condition (2) means that <<pure>>
row-wise convergence (3) takes place ``on the average'' so that that
the <<row-wise>> convergence as $k\to\infty$ is somehow coherent
with the <<principal>> convergence as $n\to\infty$.

\smallskip

{\sc Remark 1}. It can be easily verified that, since the values
under the expectation sign in (2) are nonnegative and bounded (by
two), then coherency condition (2) is equivalent to that
$$
\sup_{|t|\le T}\big|h_{n,N_n}(t)-h(t)\big|\longrightarrow 0
$$
in probability as $n\to\infty$.

\smallskip

{\sc Lemma 1.} {\it Let the family of random variables
$\{U_n\}_{n\in\mathbb{N}}$ be weakly relatively compact. Assume that
coherency condition $(2)$ holds. Then for any $t\in\mathbb{R}$ we
have}
$$
\lim_{n\to\infty}|f_n(t)-g_n(t)|=0.\eqno(4)
$$

\smallskip

{\sc Proof}. Let $\gamma\in(0,\infty)$ be a real number to be
specified later. Denote
$$
\begin{array}{rcl}
K_{1,n}&\equiv& K_{1,n}(\gamma)=\{k:\, b_{n,k}\le\gamma
d_n\},\vspace{2mm}\\
K_{2,n}&\equiv& K_{2,n}(\gamma)=\{k:\, b_{n,k}>\gamma d_n\}.
\end{array}
$$
If $t=0$, then the assertion of the lemma is trivial. Fix an
arbitrary $t\neq 0$. By the formula of total probability we have
$$
|f_n(t)-g_n(t)|=\bigg|\sum\nolimits_{k=1}^{\infty}{\sf
P}(N_n=k)\exp\Big\{it\frac{a_{n,k}-c_n}{d_n}\Big\}\Big[f_{n,k}\Big(\frac{tb_{n,k}}{d_n}\Big)-h\Big(\frac{tb_{n,k}}{d_n}\Big)\Big]\bigg|\le
$$
$$
\le\sum\nolimits_{k\in K_{1,n}}{\sf
P}(N_n=k)\Big|f_{n,k}\Big(\frac{tb_{n,k}}{d_n}\Big)-h\Big(\frac{tb_{n,k}}{d_n}\Big)\Big|+
$$
$$
+\sum\nolimits_{k\in K_{2,n}}{\sf
P}(N_n=k)\Big|f_{n,k}\Big(\frac{tb_{n,k}}{d_n}\Big)-h\Big(\frac{tb_{n,k}}{d_n}\Big)\Big|\equiv
I_1+I_2.\eqno(5)
$$
Choose an arbitrary $\epsilon>0$.

First consider $I_2$. We obviously have
$$
I_2\le\sum\nolimits_{k\in K_{2,n}(\gamma)}{\sf P}(N_n=k)={\sf
P}\big(U_n>\gamma\big).\eqno(9)
$$
The weak relative compactness of the family
$\{U_n\}_{n\in\mathbb{N}}$ implies the existence of a
$\gamma_1=\gamma_1(\epsilon)$ such that
$$
\sup_n{\sf P}\big(U_n>\gamma_1\big)<\epsilon.
$$
Therefore, setting $\gamma=\gamma_1$ from (9) we obtain
$$
I_2<\epsilon.\eqno(10)
$$
Now consider $I_1$ with $\gamma$ chosen above. If $k\in
K_{1,n}(\gamma)$, then $\big|tb_{n,k}/d_n\big|\le\gamma|t|$ and we
have
$$
I_1\le\sum\nolimits_{k\in K_{1,n}(\delta,\gamma)}{\sf
P}(N_n=k)\sup_{|\tau|\le\gamma |t|}|f_{n,k}(\tau)-h(\tau)|\le{\sf
E}\sup_{|\tau|\le\gamma |t|}|f_{n,N_n}(\tau)-h(\tau)|.
$$
Therefore, coherency condition (2) implies that there exists a
number $n_0=n_0(\epsilon,\gamma)$ such that for all $n\ge n_0$
$$
I_1<\epsilon.\eqno(11)
$$
Unifying (5), (10) and (11) we obtain that
$|f_n(t)-g_n(t)|<2\epsilon$ for $n\ge n_0$. The arbitrariness of
$\epsilon$ proves (4). The lemma is proved.

\smallskip

Lemma 1 makes it possible to use the distribution function $G_n(x)$
(see (1)) as an {\it accompanying asymptotic} approximation to
$F_n(x)\equiv{\sf P}(Z_n<x)$. In order to obtain a {\it limit}
approximation, in the next section we formulate and prove the
transfer theorem.

\section{General transfer theorem and its inversion. The structure of limit laws}

{\sc Theorem 1}. {\it Assume that coherency condition $(2)$ holds.
If there exist random variables $U$ and $V$ such that the joint
distributions of the pairs $(U_n,\,V_n)$ converge to that of the
pair $(U,\,V):$
$$
(U_n,\,V_n)\Longrightarrow (U,\,V)\ \ \ (n\to\infty),\eqno(12)
$$
then
$$
Z_n\Longrightarrow Z\eqd Y\cdot U+V\ \ \ (n\to\infty).\eqno(13)
$$
where the random variable $Y$ is independent of the pair $(U,\,V)$.}

\smallskip

{\sc Proof}. Treating $t\in\mathbb{R}$ as a fixed parameter,
represent the function $g_n(t)$ as $g_n(t)={\sf
E}h(tU_n)e^{itV_n}\equiv{\sf E}\phi_t(U_n,\,V_n)$. Since for each
$t\in\mathbb{R}$ the function $\phi_t(x,y)\equiv h(tx)e^{ity}$,
$x,y\in \mathbb{R}$, is bounded and continuous in $x$ and $y$, then
by the definition of the weak convergence we have
$$
\lim_{n\to\infty}{\sf E}\phi_t(U_n,V_n)={\sf E}\phi_t(U,V).\eqno(14)
$$
Using the Fubini theorem it can be easily verified that the function
on the right-hand side of (14) is the characteristic function of the
random variable $Y\cdot U+V$ where the random variable $Y$ is
independent of the pair $(U,V)$. Now the statement of the theorem
follows from lemma 1 by the triangle inequality. The theorem is
proved.

\smallskip

It is easy to see that relation (13) is equivalent to the following
relation between the distribution functions $F(x)$ and $H(x)$ of the
random variables $Z$ and $Y$:
$$
F(x)={\sf E}H\Big(\frac{x-V}{U}\Big),\ \ \ x\in\mathbb{R},\eqno(15)
$$
that is, the limit law for normalized randomly indexed random
variables $Z_n$ is a scale-location mixture of the distributions
which are limiting for normalized non-randomly indexed random
variables $Y_{n,k}$. Among all scale-location mixtures, {\it
variance-mean mixtures} attract a special interest (to be more
precise, we should speak of {\it normal variance-mean mixtures}).
Let us see how these mixture can appear in the double-array setting
under consideration.

Assume that the centering constants $a_{n,k}$ and $c_n$ are in some
sense proportional to the scaling constants $b_{n,k}$ and $d_n$.
Namely, assume that there exist $\rho>0$, $\alpha_n\in\mathbb{R}$
and $\beta_n\in\mathbb{R}$ such that for all $n,k\in\mathbb{N}$ we
have
$$
a_{n,k}=\frac{b_{n,k}^{\rho+1}\alpha_n}{d_n^\rho},\ \ \
c_n=d_n\beta_n,\eqno(16)
$$
and there exist finite limits
$$
\alpha=\lim_{n\to\infty}\alpha_n,\ \ \
\beta=\lim_{n\to\infty}\beta_n.
$$
Then under condition (12)
$$
(U_n,\,V_n)=\Big(\frac{b_{n,N_n}}{d_n},\,\frac{a_{n,N_n}-c_n}{d_n}\Big)=\big(U_n,\,\alpha_nU_n^{\rho+1}+\beta_n\big)\Longrightarrow
\big(U,\,\alpha U^{\rho+1}+\beta\big)\ \ \ (n\to\infty),
$$
so that in accordance with theorem 2 the limit law for $Z_n$ takes
the form
$$
{\sf P}(Z<x)={\sf E}H\Big(\frac{x-\beta-\alpha U^{\rho+1}}{U}\Big),\
\ \ x\in\mathbb{R}.
$$
If $\rho=1$, then we obtain the ``pure'' variance-mean mixture
$$
{\sf P}(Z<x)={\sf E}H\Big(\frac{x-\beta-\alpha U^2}{U}\Big),\ \ \
x\in\mathbb{R}.
$$
We will return to the discussion of convergence of randomly indexed
sequences, more precisely, of random sums, to normal scale-location
mixtures in Sect. 5.

In order to prove the result that is a partial inversion of theorem
1, for fixed random variables $Z$ and $Y$ with the characteristic
functions $f(t)$ and $h(t)$ introduce the set $\mathcal{W}(Z|Y)$
containing all pairs of random variables $(U,V)$ such that the
characteristic function $f(t)$ can be represented as
$$
f(t)={\sf E}h(tU)e^{itV},\ \ \ t\in\mathbb{R},\eqno(17)
$$
and ${\sf P}(U\ge 0)=1.$ Whatever random variables $Z$ and $Y$ are,
the set $\mathcal{W}(Z|Y)$ is always nonempty since it trivially
contains the pair $(0,Z)$. It is easy to see that representation
(17) is equivalent to relation (15) between the distribution
functions $F(x)$ and $H(x)$ of the random variables $Z$ and $Y$.

The set $\mathcal{W}(Z|Y)$ may contain more that one element. For
example, if $Y$ is the standard normal random variable and $Z\eqd
W_1-W_2$ where $W_1$ and $W_2$ are independent random variables with
the same standard exponential distribution, then along with the pair
$(0,W_1-W_2)$ the set $\mathcal{W}(Z|Y)$ contains the pair
$\big(\sqrt{W_1},0\big)$. In this case $F(x)$ is the symmetric
Laplace distribution.

Let $L_1(X_1,\,X_2)$ be the L{\'e}vy distance between the
distributions of random variables $X_1$ and $X_2$: if $F_1(x)$ and
$F_2(x)$ are the distribution functions of $X_1$ and $X_2$,
respectively, then
$$
L_1(X_1,\,X_2)=\inf\{y\ge0:\,F_2(x-y)-y\le F_1(x)\le F_2(x+y)+y\
\text{for all } x\in\mathbb{R}\}.
$$
As is well known, the L{\'e}vy distance metrizes weak convergence.
Let $L_2\big((X_1,X_2),\,(Y_1,Y_2)\big)$ be any probability metric
which metrizes weak convergence in the space of two-dimensional
random vectors. An example of such a metric is the
L{\'e}vy--Prokhorov metric (see,e. g., \cite{Zolotarev1997}).

\smallskip

{\sc Theorem 2}. {\it Let the family of random variables
$\{U_n\}_{n\in\mathbb{N}}$ be
weakly relatively compact. Assume that coherency condition $(2)$
holds. Then a random variable $Z$ such that
$$
Z_n\Longrightarrow Z\ \ \ (n\to\infty)\eqno(18)
$$
with some $c_n\in\mathbb{R}$ exists if and only if there exists a
weakly relatively compact sequence of pairs
$(U_n',\,V_n')\in\mathcal{W}(Z|Y)$, $n\in\mathbb{N}$, such that
$$
\lim_{n\to\infty}L_2\big((U_n,\,V_n),\,(U_n',\,V_n')\big)=0.\eqno(19)
$$
}

\smallskip

{\sc Proof}. {\it ``Only if'' part}. Prove that the sequence
$\{V_n\}_{n\in\mathbb{N}}$ is weakly relatively compact. The
indicator function of a set $A$ will be denoted $\mathbb{I}(A)$. By
the formula of total probability for an arbitrary $R>0$ we have
$$
{\sf P}(|V_n|>R)=\sum\nolimits_{k=1}^{\infty}{\sf
P}(N_n=k)\mathbb{I}\Big(\Big|\frac{a_{n,k}-c_n}{d_n}\Big|>R\Big)=
$$
$$
= \sum\nolimits_{k=1}^{\infty}{\sf P}(N_n=k){\sf
P}\Big(\Big|\frac{S_{n,k}-c_n}{d_n}-\frac{b_{n,k}}{d_n}\cdot\frac{S_{n,k}-a_{n,k}}{b_{n,k}}\Big|>R\Big)\le
$$
$$
\le {\sf
P}\Big(|Z_n|>\frac{R}{2}\Big)+\sum\nolimits_{k=1}^{\infty}{\sf
P}(N_n=k){\sf
P}\Big(\frac{b_{n,k}}{d_n}\cdot|Y_{n,k}|>\frac{R}{2}\Big)\equiv
I_{1,n}(R)+I_{2,n}(R).
$$
First consider $I_{2,n}(R)$. Using the set $K_{2,n}=K_{2,n}(\gamma)$
introduced in the preceding section, for an arbitrary $\gamma>0$ we
have
$$
I_{2,n}(R)=\sum\nolimits_{k\in K_{2,n}}\!\!\!{\sf P}(N_n=k){\sf
P}\Big(|Y_{n,k}|>\frac{Rd_n}{2b_{n,k}}\Big)+
\sum\nolimits_{k\notin K_{2,n}}\!\!\!{\sf P}(N_n=k){\sf
P}\Big(|Y_{n,k}|>\frac{Rd_n}{2b_{n,k}}\Big)\le
$$
$$
\le \sum\nolimits_{k\in K_{2,n}}{\sf P}(N_n=k){\sf
P}\Big(|Y_{n,k}|>\frac{R}{2\gamma}\Big)+{\sf
P}\big(U_n>\gamma\big)\le {\sf
P}\Big(|Y_{n,N_n}|>\frac{R}{2\gamma}\Big)+{\sf
P}\big(U_n>\gamma\big) .\eqno(20)
$$
Fix an arbitrary $\epsilon>0$. Choose $\gamma=\gamma(\epsilon)$ so
that
$$
{\sf P}\big(U_n>\gamma(\epsilon)\big)<\epsilon\eqno(21)
$$
for all $n\in\mathbb{N}$. This is possible due to the weak relative
compactness of the family $\{U_n\}_{n\in\mathbb{N}}$. Now choose
$R'=R'(\epsilon)$ so that
$$
{\sf
P}\Big(|Y_{n,N_n}|>\frac{R'(\epsilon)}{2\gamma(\epsilon)}\Big)<\epsilon.\eqno(22)
$$
This is possible due to the weak relative compactness of the family
$\{Y_{n,N_n}\}_{n\in\mathbb{N}}$ implied by coherency condition (2).
Thus, from (20), (21) and (22) we obtain
$$
I_{2,n}\big(R'(\epsilon)\big)<2\epsilon\eqno(23)
$$
for all $n\in\mathbb{N}$. Now consider $I_{1,n}(R)$. From (18) it
follows that there exists an $R''=R''(\epsilon)$ such that
$$
I_{1,n}\big(R''(\epsilon)\big)<\epsilon\eqno(24)
$$
for all $n\in\mathbb{N}$. From (23) and (24) it follows that if
$R>\max\big\{R',\,R''\big\}$, then
$$
\sup_{n}{\sf P}(|V_n|>R)<3\epsilon
$$
and by virtue of the arbitrariness of $\epsilon>0$, the family
$\{V_n\}_{n\in\mathbb{N}}$ is weakly relatively compact. Hence, the
family of pairs $\{(U_n,\,V_n)\}_{n\in\mathbb{N}}$ is weakly
relatively compact.

Denote
$$
\epsilon_n=\inf\big\{L_2\big((U_n,\,V_n),\,(U,\,V)\big):\
(U,\,V)\in\mathcal{W}(Z|Y)\big\},\ \  n=1,2,...
$$
Prove that $\epsilon_n\to 0$ as $n\to\infty$. Assume the contrary.
In this case $\epsilon_n\ge M$ for some $M>0$ and all $n$ from some
subsequence $\mathcal{N}$ of natural numbers. Choose a subsequence
$\mathcal{N}_1\subseteq\mathcal{N}$ so that the sequence of pairs
$\{(U_n,\,V_n)\}_{n\in\mathcal{N}_1}$ weakly converges to some pair
$(U,\,V)$. As this is so, for all $n\in\mathcal{N}_1$ large enough
we will have $L_2\big((U_n,\,V_n),\,(U,V)\big)<M$. Applying theorem
1 to the sequence $\{(U_n,\,V_n)\}_{n\in\mathcal{N}_1}$ we make sure
that $(U,V)\in \mathcal{W}(Z|Y)$ since condition (18) implies the
coincidence of the limits of all convergent subsequences of
$\{Z_n\}$. We arrive at the contradiction with the assumption that
$\epsilon_n>M$ for all $n\in\mathcal{N}_1$. Hence, $\epsilon_n\to 0$
as $n\to\infty$. For each $n\in\mathbb{N}$ choose a pair
$(U_n',\,V_n')\in\mathcal{W}(Z|Y)$ such that
$$
L_2\big((U_n,\,V_n),\,(U_n',\,V_n')\big)\le\epsilon_n+{\textstyle\frac1n}.
$$
The sequence $\{(U_n',\,V_n')\}_{n\in\mathbb{N}}$ obviously
satisfies condition (19). Its weak relative compactness follows from
(19) and the weak relative compactness of the sequence
$\{(U_n,\,V_n)\}_{n\in\mathbb{N}}$ established above.

\smallskip

{\it ``If'' part}. Assume that the sequence
$\{Z_n\}_{n\in\mathbb{N}}$ does not converge weakly to $Z$ as
$n\to\infty$. In that case the inequality $L_1(Z_n,\,Z)\ge M$ holds
for some $M>0$ and all $n$ from some subsequence $\mathcal{N}$ of
natural numbers. Choose a subsequence
$\mathcal{N}_1\subseteq\mathcal{N}$ so that the sequence of pairs
$\{(U_n',\,V_n')\}_{n\in\mathcal{N}_1}$ weakly converges to some
pair $(U,\,V)$. Repeating the reasoning used to prove theorem 1 we
make sure that for any $t\in\mathbb{R}$
$$
{\sf E}e^{itZ}={\sf E}h(tU_n')e^{itV_n'}\longrightarrow{\sf
E}h(tU)e^{itV}
$$
as $n\to\infty$, $n\in\mathcal{N}_1$, that is,
$(U,V)\in\mathcal{W}(Z|Y)$. From the triangle inequality
$$
L_2\big((U_n,\,V_n),\,(U,V)\big)\le
L_2\big((U_n,\,V_n),\,(U_n',\,V_n')\big)+L_2\big((U_n',\,V_n'),\,(U,V)\big)
$$
and condition (19) it follows that
$L_2\big((U_n,\,V_n),\,(U,V)\big)\to0$ as $n\to\infty$,
$n\in\mathcal{N}_1$. Apply theorem 1 to the double array
$\{Y_{n,k}\}_{k\in\mathbb{N},\,n\in\mathcal{N}_1}$ and the sequence
$\{(U_n,\,V_n)\}_{n\in\mathcal{N}_1}$. As a result we obtain that
$L_1(Z_n,Z)\to0$ as $n\to\infty$, $n\in\mathcal{N}_1$, contradicting
the assumption that $L_1(Z_n,Z)\ge M>0$ for $n\in\mathcal{N}_1$.
Thus, the theorem is completely proved.

\smallskip

{\sc Remark 2}. It should be noted that in \cite{Korolev1993} and
some subsequent papers a stronger and less convenient version of the
coherency condition was used. Furthermore, in \cite{Korolev1993} and
the subsequent papers the statements analogous to lemma 1 and
theorems 1 and 2 were proved under the additional assumption of the
weak relative compactness of the family
$\{Y_{n,k}\}_{n,k\in\mathbb{N}}$.

\section{Limit theorems for random sums of independent random variables}

Let $\{X_{n,j}\}_{j\ge1}$, $n\in\mathbb{N},$ be a double array of
row-wise independent not necessarily identically distributed random
variables. For $n,k\in\mathbb{N}$ denote
$$
S_{n,k}=X_{n,1}+...+X_{n,k}.\eqno(25)
$$
If $S_{n,k}$ is a sum of independent random variables, then the
condition of weak relative compactness of the sequence
$\{U_n\}_{n\in\mathbb{N}}$ used in the preceding section can be
replaced by the condition of weak relative compactness of the family
$\{Y_{n,k}\}_{n,k\in\mathbb{N}}$ which is in fact considerably less
restrictive. Indeed, let, for example, the random variables
$S_{n,k}$ possess moments of some order $\delta>0$. Then, if we
choose $b_{n,k}=({\sf E}|S_{n,k}-a_{n,k}|^{\delta})^{1/\delta}$,
then by the Markov inequality
$$
\lim_{R\to\infty}\sup_{n,k\in\mathbb{N}}{\sf
P}\big(|Y_{n,k}|>R\big)\le\lim_{R\to\infty}\frac{1}{R\,^{\delta}}=0,
$$
that is, the family $\{Y_{n,k}\}_{n,k\in\mathbb{N}}$ is weakly
relatively compact.

\smallskip

{\sc Theorem 3}. {\it Assume that the random variables $S_{n,k}$
have the form $(25)$. Let the family of random variables
$\{Y_{n,k}\}_{n,k\in\mathbb{N}}$ be weakly relatively compact.
Assume that coherency condition $(2)$ holds. Then convergence $(18)$
of normalized random sums $Z_n$ to some random variable $Z$ takes
place with some $c_n\in\mathbb{R}$ if and only if there exists a
weakly relatively compact sequence of pairs
$(U_n',\,V_n')\in\mathcal{W}(Z|Y)$, $n\in\mathbb{N}$, such that
condition $(19)$ holds. }

\smallskip

{\sc Proof}. It suffices to prove that in the case under
consideration condition (18) implies the weak relative compactness
of the family $\{U_n\}_{n\in\mathbb{N}}$. In what follows the
symmetrization of a random variable $X$ will be denoted $X^{(s)}$,
$X^{(s)}= X-X'$ where $X'$ is a random variable independent of $X$
such that $X'\eqd X$. For $q\in(0,1)$ let $\ell_n(q)$ be the
greatest lower bound of $q$-quantiles of the random variable $N_n$,
$n\in\mathbb{N}$. Assume that for each $n$ the random variables
$N_n,X_{n,1}^{(s)},X_{n,2}^{(s)},...$ are jointly independent and
introduce the random variables
$$
Q_n=\frac{1}{d_n}\sum\nolimits_{j=1}^{N_n}X_{n,j}^{(s)},\ \ \
n\in\mathbb{N}.
$$
Using the symmetrization inequality
$$
{\sf P}\big(X^{(s)}\ge
R\big)\le 2\,{\sf P}(|X-a|\ge\textstyle{\frac{R}{2}})
$$
which is valid for any random variable $X$, any $a\in\mathbb{R}$ and
$R>0$ (see, e. g., \cite{Loeve1977}), we obtain
$$
{\sf P}(|Q_n|\ge R)=\sum\nolimits_{k=1}^{\infty}{\sf P}(N_n=k){\sf
P}\Big(\Big|\frac{1}{d_n}\sum\nolimits_{j=1}^kX_{n,j}^{(s)}\Big|\ge
R\Big)\le
$$
$$
\le 2\sum\nolimits_{k=1}^{\infty}{\sf P}(N_n=k){\sf
P}\Big(\Big|\frac{1}{d_n}\Big(\sum\nolimits_{j=1}^kX_{n,j}-c_n\Big)\Big|\ge
\frac{R}{2}\Big)=2\,{\sf
P}\Big(\Big|\frac{1}{d_n}\Big(\sum\nolimits_{j=1}^{N_n}X_{n,j}-c_n\Big)\Big|\ge
\frac{R}{2}\Big)
$$
for any $R>0$ and $n\in\mathbb{N}$. Hence,
$$
\lim_{R\to\infty}\sup_n{\sf P}(|Q_n|\ge R)\le
2\lim_{R\to\infty}\sup_n{\sf
P}\Big(\Big|\frac{1}{d_n}\Big(\sum\nolimits_{j=1}^{N_n}X_{n,j}-c_n\Big)\Big|\ge
\frac{R}{2}\Big)=0
$$
by virtue of (18). Hence, the sequence $\{Q_n\}_{n\in\mathbb{N}}$ is
weakly relatively compact.


Now prove that
$$
C(q)\equiv\sup_n\frac{b_{n,\ell_n(q)}}{d_n}<\infty\eqno(26)
$$
for each $q\in(0,1)$. For this purpose we use the L{\'e}vy
inequality
$$
{\sf P}\Big(\max_{1\le m\le
k}\Big|\sum\nolimits_{j=1}^mX_j^{(s)}\Big|\ge R\Big)\le2\,{\sf
P}\Big(\Big|\sum\nolimits_{j=1}^kX_j^{(s)}\Big|\ge R\Big)
$$
which is valid for any independent random variables $X_1,...,X_k$
and any $R>0$ and for an arbitrary $q\in(0,1)$ obtain the following
chain of inequalities:
$$
2\,{\sf P}(|Q_n|\ge R)=2\sum\nolimits_{k=1}^{\infty}{\sf
P}(N_n=k){\sf
P}\Big(\Big|\frac{1}{d_n}\sum\nolimits_{j=1}^kX_{n,j}^{(s)}\Big|\ge
R\Big)\ge
$$
$$
\ge 2\sum\nolimits_{k\ge\ell_n(q)}{\sf P}(N_n=k){\sf
P}\Big(\Big|\frac{1}{d_n}\sum\nolimits_{j=1}^kX_{n,j}^{(s)}\Big|\ge
R\Big)\ge
$$
$$
\ge \sum\nolimits_{k\ge\ell_n(q)}\!\!{\sf P}(N_n=k){\sf
P}\Big(\Big|\frac{1}{d_n}\sum\nolimits_{j=1}^{\ell_n(q)}X_{n,j}^{(s)}\Big|\ge
R\Big)\!\!={\sf P}\big(N_n\ge\ell_n(q)\big){\sf
P}\Big(\Big|\frac{1}{d_n}\sum\nolimits_{j=1}^{\ell_n(q)}X_{n,j}^{(s)}\Big|\ge
R\Big)=
$$
$$
=(1-q){\sf
P}\Big(\Big|\frac{1}{d_n}\sum\nolimits_{j=1}^{\ell_n(q)}X_{n,j}^{(s)}\Big|\ge
R\Big).
$$
Hence, the weak relative compactness of the family
$\{Q_n\}_{n\in\mathbb{N}}$ established above, for each $q\in(0,1)$
implies the weak relative compactness of the family
$\{Q_n^{(q)}\}_{n\in\mathbb{N}}$ where
$$
Q_n^{(q)}=\frac{1}{d_n}\sum\nolimits_{j=1}^{\ell_n(q)}X_{n,j}^{(s)},\
\ \ n\in\mathbb{N}.
$$
Assume that (26) does not hold. In that case there exist a
$q^*\in(0,1)$ and a sequence $\mathcal{N}$ of natural numbers such
that
$$
\frac{b_{n,\ell_n(q^*)}}{d_n}\longrightarrow \infty,\ \ \
n\to\infty,\ n\in\mathcal{N}.\eqno(27)
$$
According to the conditions of the theorem, the family of random
variables
$\big\{Y_{n,k}=(S_{n,k}-a_{n,k})/b_{n,k}\big\}_{n,k\in\mathbb{N}}$
is weakly relatively compact. Therefore, a subsequence
$\mathcal{N}_1\subseteq\mathcal{N}$ can be chosen so that
$$
Y_{n,\ell_n(q^*)}=\frac{1}{b_{n,\ell_n(q^*)}}\Big(\sum\nolimits_{j=1}^{\ell_n(q^*)}X_{n,j}-a_{n,\ell_n(q^*)}\Big)\Longrightarrow
Y,\ \ \ n\to\infty,\ n\in\mathcal{N}_1,\eqno(28)
$$
where $Y$ is some random variable. From (27) and (28) it follows
that for any $R\in\mathbb{R}$
$$
{\sf P}(Q_n^{(q^*)}\le R)={\sf
P}\bigg(Y_{n,\ell_n(q^*)}^{(s)}\le\frac{d_nR}{b_{n,\ell_n(q^*)}}\bigg)\longrightarrow
{\sf P}(Y^{(s)}\le 0)\ge\textstyle{\frac12},\ \ \ n\to\infty,\
n\in\mathcal{N}_1,
$$
contradicting the weak relative compactness of the family
$\{Q_n^{(q^*)}\}$ established above. So, (26) holds for any
$q\in(0,1)$.

It is easy to make sure that $N_n\eqd\ell_n(W)$ where $W$ is a
random variable with the uniform distribution on $(0,1)$. Therefore,
with the account of (26) for any $R\ge0$ and $n\in\mathbb{N}$ we
have
$$
{\sf P}(U_n\ge R)={\sf P}\Big(\frac{b_{n,\ell_n(W)}}{d_n}\ge
R\Big)=\int_{0}^{1}\mathbb{I}\Big(\frac{b_{n,\ell_n(q)}}{d_n}\ge
R\Big)dq\le
$$
$$
\le\int_{0}^{1}\mathbb{I}\big(C(q)\ge R\big)dq={\sf P}\big(C(W)\ge
R\big)
$$
so that
$$
\lim_{R\to\infty}\sup_n{\sf P}\big(U_n\ge
R\big)=\lim_{R\to\infty}{\sf P}\big(C(W)\ge R\big)=0,
$$
that is, the sequence $\{U_n\}_{n\in\mathbb{N}}$ is weakly
relatively compact.

The rest of the proof of theorem 3 repeats that of theorem 2
word-for-word. The theorem is proved.

\section{A version of the central limit theorem for random sums
with a normal variance-mean mixture as the limiting law}

Let $\{X_{n,j}\}_{j\ge1}$, $n\in\mathbb{N},$ be a double array of
row-wise independent not necessarily identically distributed random
variables. As in the preceding section, let
$$
S_{n,k}=X_{n,1}+...+X_{n,k},\ \ \ n,k\in\mathbb{N}.
$$
The distribution function and the characteristic function of the
random variable $X_{n,j}$ will be denoted $F_{n,j}(x)$ and
$f_{n,j}(t)$, respectively,
$$
f_{n,j}(t)=\int_{-\infty}^{\infty}e^{itx}dF_{n,j}(x),\ \ \
t\in\mathbb{R}.
$$
It is easy to see that in this case
$$
h_{n,k}(t)\equiv {\sf
E}\exp\big\{itY_{n,k}\big\}=\exp\Big\{-it\frac{a_{n,k}}{b_{n,k}}\Big\}\prod\nolimits_{j=1}^kf_{n,j}\Big(\frac{t}{b_{n,k}}\Big),\
\ t\in\mathbb{R}.
$$

Denote $\mu_{n,j}={\sf E}X_{n,j}$, $\sigma_{n,j}^2={\sf D}X_{n,j}$
and assume that $0<\sigma_{n,j}^2<\infty$, $n,j\in\mathbb{N}$.
Denote
$$
A_{n,k}=\mu_{n,1}+...+\mu_{n,k}\ \big(={\sf E}S_{n,k}\big),\ \ \
B_{n,k}^2=\sigma_{n,1}^2+...+\sigma_{n,k}^2\ \big(={\sf
D}S_{n,k}\big)
$$
It is easy to make sure that ${\sf E}S_{n,N_n}={\sf E}A_{n,N_n}$,
${\sf D}S_{n,N_n}={\sf E}B_{n,N_n}^2+{\sf D}A_{n,N_n}$,
$n\in\mathbb{N}$. In order to formulate a version of the central
limit theorem for random sums with the limiting distribution being a
normal variance-mean mixture, assume that non-random sums, as usual,
are centered by their expectations and normalized by by their mean
square deviations and put $a_{n,k}=A_{n,k}$,
$b_{n,k}=\sqrt{B_{n,k}^2}$, $n,k\in\mathbb{N}$. Although it would
have been quite natural to normalize random sums by their mean
square deviations as well, for simplicity we will use slightly
different normalizing constants and put $d_n=\sqrt{{\sf
E}B_{n,N_n}^2}$.

Let $\Phi(x)$ be the standard normal distribution function,
$$
\Phi(x)=\frac{1}{\sqrt{2\pi}}\int_{-\infty}^{x}e^{-z^2/2}dz,\ \ \
x\in\mathbb{R}.
$$

\smallskip

{\sc Theorem 4.} {\it Assume that the following conditions
hold{\rm:}

\noindent $(i)$ for every $n\in\mathbb{N}$ there exist real numbers
$\alpha_n$ such that
$$
\mu_{n,j}=\frac{\alpha_n\sigma_{n,j}^2}{\sqrt{{\sf E}B_{n,N_n}^2}},\
\ \ n,j\in\mathbb{N},\eqno(29)
$$
and
$$
\lim_{n\to\infty}\alpha_n=\alpha,\ \ \ 0<|\alpha|<\infty;\eqno(30)
$$

\noindent $(ii)$ $($the random Lindeberg condition$)$ for any
$\epsilon>0$
$$
\lim_{n\to\infty}{\sf
E}\frac{1}{B_{n,N_n}^2}\sum\nolimits_{j=1}^{N_n}\int_{|x-\mu_{n,j}|>\epsilon
B_{n,N_n}}\!\!\!(x-\mu_{n,j})^2dF_{n,j}(x)=0.\eqno(31)
$$

\noindent Then the convergence of the normalized random sums
$$
\frac{S_{n,N_n}}{\sqrt{{\sf E}B_{n,N_n}^2}}\Longrightarrow
Z\eqno(32)
$$
to some random variable $Z$ as $n\to\infty$ takes place if and only
if there exists a random variable $U$ such that
$$
{\sf P}(Z<x)={\sf E}\Phi\Big(\frac{x-\alpha U}{\sqrt{U}}\Big),\ \
x\in\mathbb{R},\eqno(33)
$$
and
$$
\frac{B_{n,N_n}^2}{{\sf E}B_{n,N_n}^2}\Longrightarrow U\ \ \
(n\to\infty).\eqno(34)
$$
}

\smallskip

{\sc Proof}. We will deduce theorem 4 as a corollary of theorem 3.

First, let $a_{n,k}=A_{n,k}$, $b_{n,k}=B_{n,k}$, $n,k\in\mathbb{N}$.
Then the family of the random variables
$\{Y_{n,k}\}_{n,k\in\mathbb{N}}$ is weakly relatively compact, since
by the Chebyshev inequality
$$
\lim_{R\to\infty}\sup_{n,k\in\mathbb{N}}{\sf
P}\big(|Y_{n,k}|>R\big)=\lim_{R\to\infty}\sup_{n,k\in\mathbb{N}}{\sf
P}\bigg(\bigg|\frac{S_{n,k}-A_{n,k}}{B_{n,k}}\bigg|>R\bigg)\le\lim_{R\to\infty}\frac{1}{R^2}=0.
$$

Second, prove that under the conditions of the theorem the coherency
condition (2) holds with $h(t)=e^{-t^2/2}$, $t\in\mathbb{R}$. Denote
$\Delta_{n,k}(x)=|H_{n,k}(x)-\Phi(x)|$ where $H_{n,k}(x)={\sf
P}(S_{n,k}-A_{n,k}<B_{n,k}x)$. By integration by parts, for any
$t\in\mathbb{R}$ we have
$$
\big|h_{n,k}(t)-e^{-t^2/2}\big|=\bigg|it\int_{-\infty}^{\infty}e^{itx}\big[H_{n,k}(x)-\Phi(x)\big]dx\bigg|\le
|t|\int_{-\infty}^{\infty}\Delta_{n,k}(x)dx,\eqno(35)
$$
To estimate the integrand on the right-hand side of (35) we will use
the following result of V. V. Petrov \cite{Petrov1979}. Let
$\mathcal{G}$ be the class of real-valued functions $g(x)$ of the
argument $x\in\mathbb{R}$ such that the function $g(x)$ is even,
nonnegative for all $x$ and positive for $x>0$; the functions $g(x)$
and $x/g(x)$ are nondecreasing for $x>0$. In \cite{Petrov1979} it
was proved that, whatever a function $g\in\mathcal{G}$ is, if ${\sf
E} X_{n,j}^2g(X_{n,j})<\infty$, $n,j\in\mathbb{N}$, then there
exists a positive finite absolute constant $C$ such that for any
$x\in\mathbb{R}$
$$
\Delta_{n,k}(x)\le\frac{C}{B_{n,k}^2(1+|x|)^2g\big(B_{n,k}(1+|x|)\big)}\sum\nolimits_{j=1}^k{\sf
E}(X_{n,j}-\mu_{n,j})^2g(X_{n,j}-\mu_{n,j}).
$$
Hence, it is easy to see that the properties of the function
$g\in\mathcal{G}$ guarantee that
$$
\Delta_{n,k}(x)\le\frac{C}{(1+|x|)^2}\cdot\frac{1}{B_{n,k}^2g\big(B_{n,k})\big)}\sum\nolimits_{j=1}^k{\sf
E}(X_{n,j}-\mu_{n,j})^2g(X_{n,j}-\mu_{n,j}).\eqno(36)
$$
Now choosing $g(x)=\min\{|x|,\,B_{n,k}\}\in\mathcal{G}$ and
repeating the reasoning used to prove theorem 7 in Sect. 3, Chapt. V
of \cite{Petrov1987} and relation (3.8) there, from (36) we obtain
that for all $n,k\in\mathbb{N}$ and an arbitrary $\epsilon>0$
$$
\Delta_{n,k}(x)\le
\frac{2C}{(1+|x|)^2}\bigg\{\epsilon+\frac{1}{B_{n,k}^2}
\sum\nolimits_{j=1}^{k}{\sf
E}\big[(X_{n,j}-\mu_{n,j})^2\mathbb{I}(|X_{n,j}-\mu_{n,j}|>\epsilon
B_{n,k})\big]\bigg\}.\eqno(37)
$$
Using (35) and (37) we obtain that for arbitrary $\epsilon>0$ and
$T\in(0,\infty)$
$$
{\sf E}\sup_{|t|\le T}|h_{n,N_n}(t)-e^{-t^2/2}|\le |T|{\sf
E}\int_{-\infty}^{\infty}\big|H_{n,N_n}(x)-\Phi(x)\big|dx\le
$$
$$
\le 4C|T|\bigg\{\epsilon+{\sf E}\frac{1}{B_{n,N_n}^2}
\sum\nolimits_{j=1}^{N_n} \int_{|x-\mu_{n,j}|>\epsilon
B_{n,N_n}}\!\!\!(x-\mu_{n,j})^2dF_{n,j}(x)\bigg\}.
$$
Hence, from (31) it follows that for an arbitrary $\epsilon>0$
$$
\lim_{n\to\infty}{\sf E}\sup_{|t|\le T}|h_{n,N_n}(t)-e^{-t^2/2}|\le
4C|T|\epsilon,
$$
and since $\epsilon>0$ can be taken arbitrarily small, the coherency
condition (2) holds.

Third, let $d_n=\sqrt{{\sf E}B_{n,N_n}^2}$, $c_n=0$,
$n\in\mathbb{N}$. As this is so, relations (29) and (30) guarantee
that relation (16) holds with $\rho=1$, $\beta_n=\beta=0$,
$n\in\mathbb{N}$, so that if (34) holds along some subsequence
$\mathcal{N}$ of natural numbers, then
$$
(U_n,\,V_n)=\bigg(\sqrt{\frac{B_{n,N_n}^2}{{\sf
E}B_{n,N_n}^2}},\,\alpha_n\frac{B_{n,N_n}^2}{{\sf
E}B_{n,N_n}^2}\bigg)\Longrightarrow \big(\sqrt{U},\,\alpha U\big),\
\ \ n\to\infty,\ n\in\mathcal{N},
$$
so that the limit law has the form of normal variance-mean mixture
(33).

Fourth, recently in \cite{Korolev2013} it was proved that normal
variance-mean mixtures are identifiable, that is, if ${\sf
P}(Y<x)=\Phi(x)$, then the set $\mathcal{W}(Z|Y)$ contains at most
one pair of the form $\big(\sqrt{U},\,\alpha U\big)$. This means
that in the case under consideration condition (19) reduces to (34).
The theorem is proved.

\smallskip

{\sc Remark 3}. In accordance with what has been said in remark 1,
the random Lindeberg condition $(ii)$ can be used in the following
form: for any $\epsilon>0$
$$
\frac{1}{B_{n,N_n}^2}\sum\nolimits_{j=1}^{N_n}\int_{|x-\mu_{n,j}|>\epsilon
B_{n,N_n}}\!\!\!(x-\mu_{n,j})^2dF_{n,j}(x)\longrightarrow 0
$$
in probability as $n\to\infty$.

\smallskip

{\sc Remark 4}. For $n,j,k\in\mathbb{N}$ denote $\nu_{n,j}^3={\sf
E}|X_{n,j}-\mu_{n,j}|^3$,
$M_{n,k}^3=\nu_{n,1}^3+\ldots+\nu_{n,k}^3$,
$L_{n,k}^3=M_{n,k}^3B_{n,k}^{-3}$. It is easy to see that for each
$n\in\mathbb{N}$
$$
{\sf
E}\frac{1}{B_{n,N_n}^2}\sum\nolimits_{j=1}^{N_n}\int_{|x-\mu_{n,j}|>\epsilon
B_{n,N_n}}\!\!\!(x-\mu_{n,j})^2dF_{n,j}(x)\le
$$
$$
\le\sum\nolimits_{k=1}^{\infty}{\sf P}(N_n=k)\frac{1}{\epsilon
B_{n,k}^3}\sum\nolimits_{j=1}^{k}\int_{|x-\mu_{n,j}|>\epsilon
B_{n,k}}\!\!\!|x-\mu_{n,j}|^3dF_{n,j}(x)\le
$$
$$
\le\frac{1}{\epsilon}\sum\nolimits_{k=1}^{\infty}{\sf
P}(N_n=k)\frac{M_{n,k}^3}{B_{n,k}^3}=\frac{1}{\epsilon}{\sf
E}L_{n,N_n}^3.
$$
Therefore, if the third absolute moments of the summands exist, then
the random Lindeberg condition $(ii)$ follows from {\it the random
Lyapunov condition}
$$
\lim_{n\to\infty}{\sf E}L_{n,N_n}^3=0
$$
which seems to be more easily verifiable than the random Lindeberg
condition. For example, let for each $n\in\mathbb{N}$ the random
variables $X_{n,1},X_{n,2},\ldots$ be identically distributed.
Denote the generating function of the random variable $N_n$ by
$\psi_n(s)$,
$$
\psi_n(s)={\sf E}s^{N_n}=\sum\nolimits_{k=1}^{\infty}s^k{\sf
P}(N_n=k),\ \ 0\le s\le1.
$$
Then
$$
{\sf E}L_{n,N_n}=\frac{\nu_{n,1}^3}{\sigma_{n,1}^3}{\sf
E}\frac{1}{\sqrt{N_n}}\le\frac{\nu_{n,1}^3}{\sigma_{n,1}^3}\sqrt{{\sf
E}\frac{1}{N_n}}=\frac{\nu_{n,1}^3}{\sigma_{n,1}^3}\Big(\int_0^1\frac{\psi_n(s)}{s}ds\Big)^{1/2},
$$
and the random Lyapunov condition holds, if
$$
\lim_{n\to\infty}\frac{\nu_{n,1}^3}{\sigma_{n,1}^3}\Big(\int_0^1\frac{\psi_n(s)}{s}ds\Big)^{1/2}=0
$$
implying the random Lindeberg condition and hence, the coherency
condition.

\smallskip

The class of normal variance-mean mixtures is very wide and, in
particular, contains the class of generalized hyperbolic
distributions which, in turn, contains (a)~symmetric and skew
Student distributions (including the Cauchy distribution) with
inverse gamma mixing distributions; (b)~variance gamma distributions
(including symmetric and non-symmetric Laplace distributions) with
gamma mixing distributions;
(c)~normal$\backslash\!\backslash$inverse Gaussian distributions
with inverse Gaussian mixing distributions including symmetric
stable laws. By variance-mean mixing many other initially symmetric
types represented as pure scale mixtures of normal laws can be
skewed, e. g., as it was done to obtain non-symmetric exponential
power distributions in \cite{GK2013}.

According to theorem 4, all these laws can be limiting for random
sums of independent non-identically distributed random variables.
For example, to obtain the skew Student distribution for $Z$ it is
necessary and sufficient that in (33) and (34) the random variable
$U$ has the inverse gamma distribution \cite{KorolevSokolov2012}. To
obtain the variance gamma distribution for $Z$ it is necessary and
sufficient that in (33) and (34) the random variable $U$ has the
gamma distribution \cite{KorolevSokolov2012}. In particular, for $Z$
to have the asymmetric Laplace distribution it is necessary and
sufficient that $U$ has the exponential distribution.

\smallskip

{\sc Remark 5}. Note that the non-random sums in the coherency
condition are centered, whereas in (32) the random sums are not
centered, and if $\alpha\neq 0$, then the limit distribution for
random sums becomes skew unlike usual non-random summation, where
the presence of the systematic bias of the summands results in that
the limit distribution becomes just shifted. So, if non-centered
random sums are used as models of some real phenomena and the limit
variance-mean mixture is skew, then it can be suspected that the
summands are actually biased.

\smallskip

{\sc Remark 6}. In limit theorems of probability theory and
mathematical statistics, centering and normalization of random
variables are used to obtain non-trivial asymptotic distributions.
It should be especially noted that to obtain reasonable
approximation to the distribution of the basic random variables (in
our case, $S_{n,N_n}$), both centering and normalizing values should
be non-random. Otherwise the approximate distribution becomes random
itself and, say, the problem of evaluation of quantiles becomes
senseless.

\renewcommand{\refname}{References}

\end{document}